\documentclass[12pt]{amsart}
\usepackage[a4paper]{geometry}
 \newtheorem{thm}{Theorem}[section]
 
 \newtheorem{lem}[thm]{Lemma}
 
 \theoremstyle{definition}
 
 \theoremstyle{remark}
 
 \numberwithin{equation}{section}

\begin{document}
\title {\textbf{Biconservative ideal hypersurfaces in Euclidean spaces}}
\author{Deepika, Andreas Arvanitoyeorgos}

 \maketitle
\begin{abstract}
A biconservative submanifold of a Riemannian manifold is a sub-
manifold with divergence free stress-energy tensor with respect to
bienergy. These are generalizations of biharamonic submanifolds. In
2013, B. Y. Chen and M.I. Munteanu proved that $\delta(2)$-ideal and
$\delta(3)$-ideal biharmonic hypersurfaces in Euclidean space are
minimal. In this paper, we generalize this result for
$\delta(2)$-ideal and $\delta(3)$-ideal bisonservative hypersurfaces
in Euclidean space. Also, we study $\delta(4)$-ideal biconservative
hypersurfaces in Euclidean space $\mathbb{E}^{6}$ having constant
scalar curvature. We prove that such a hypersurface must be of
constant mean curvature.
\end{abstract}

\maketitle
\section{\textbf{Introduction}}

Recently, the theory of biconservative submanifolds, which is
closely related to biharmonic submanifolds, is an active area of
research in differential geometry. A biharmonic map $\varphi:
(M,g)\rightarrow (N,h)$ is a critical point of the bienergy
functional
$E_{2}(\varphi)=\frac{1}{2}\int_{M}|\tau(\varphi)|^{2}v_{g}$, where
$\tau(\varphi)$ is the tension field of $\varphi$. These critical
points are given by the vanishing of the bitension field,
\emph{i.e.}
\begin{center}
$\tau_{2}(\varphi)= -\Delta \tau (\varphi)- \mbox{trace} R^{N}(d
\varphi, \tau (\varphi))d \varphi =0,$
\end{center}
where $R^{N}$ is the curvature tensor of $N$.

As described in \cite{SCA}, the stress energy tensor for bienergy is
defined as
\begin{center}$\begin{array}{rcl}
S_{2}(X,Y)=\frac{1}{2}|\tau(\varphi)|^{2}\langle X,Y\rangle+\langle
d\varphi, \nabla\tau(\varphi)\rangle \langle X,Y\rangle-\langle
d\varphi(X), \nabla_{Y}\tau(\varphi)\rangle \\-\langle d\varphi(Y),
\nabla_{X}\tau(\varphi)\rangle\end{array}$
\end{center}
and it satisfies
\begin{center}
$\mbox{div } S_{2}= -\langle \tau_{2}(\varphi), d\varphi\rangle,$
\end{center}
thus conforming to the principle of a stress-energy tensor for the
bienergy. If $\varphi$ is an isometric immersion with div $S_{2}=0$
then tangent part of the corresponding bitension field vanishes.

The concept of biconservative comes from the conservativity of the
stress-energy tensor $S_{2}$ for bienergy, $i.e.$ $\mbox{div
}S_{2}=0$. In fact, we can say that isometric immersion $\varphi:M
\rightarrow N$ is called biconservative if the tangential part of
bitension field vanishes. Thus, biharmonicity always implies
biconservativity.

It can be easily seen that a biconservative hypersurface $M^{n}$ in
a Riemannian manifold $N^{n+1}$ satisfies (\cite{BY1}, \cite{RSC})

\begin{equation}\label{1.1}
 2\mathcal{A} (\rm grad \emph{H})+ \emph{nH} \hspace{.1 cm} grad \emph{H} = 2 \emph{H}\hspace{.1 cm}
 Ricci^{N}(\xi)^{\top},
\end{equation}
where $\mathcal{A}$ is the shape operator, $H$ is the mean curvature
function and $\mbox{Ricci}^{N}(\xi)^{\top}$ is the tangent component
of the Ricci curvature of $N$ in the direction of the unit normal
$\xi$ of $M^{n}$ in $N^{n+1}$.

In this paper, we consider a biconservative hypersurface $M^{n}$ in
the Euclidean space $\mathbb{E}^{n+1}$. In this case (\ref{1.1})
becomes

\begin{equation}\label{1.2}
 2\mathcal{A} (\rm grad \emph{H})+ \emph{nH} \hspace{.1 cm} grad \emph{H} =
 0,
\end{equation}
which is the tangential component of $\triangle \vec {H} = 0$, where
$\triangle$ is a Laplace operator. This paper will help to study
much larger family of hypersurfaces including biharmonic
hypersurfaces in Euclidean space.

From (\ref{1.2}), it is obvious that hypersurfaces with constant
mean curvature are always biconservative. The question that arises
is whether there exist biconservative hypersurfaces which are not of
constant mean curvature, known as proper biconservative.

The concept of biconservative hypersurfaces have been studied by
several geometers. The first result on biconservative hypersurfaces
was obtained by T. Hasanis and T. Vlachos in \cite{TT}, who called
them as H-hypersurfaces. In \cite{RSC} R. Caddeo et al. introduced
the notion of biconservative and proved that a biconservative
surface in Euclidean 3-space is either a surface of constant mean
curvature or a surface of revolution (cf. \cite{TT, YF}). In
\cite{BY3} the authors proved that a $\delta(2)-$ideal
biconservative hypersurface in Euclidean space $\mathbb{E}^{n}$ $(n
\geq 3)$ (see definition below) is either minimal or a spherical
hypercylinder. In \cite{SCA} Montaldo et al. studied proper
$SO(p+1)\times SO(p+1)$-invariant biconservative hypersurfaces and
proper $SO(p+1)$-invariant biconservative hypersurfaces in Euclidean
space $\mathbb{E}^{n}$. Also, Fectu et al. classified biconservative
surfaces in $\mathbb{S}^{n}\times \mathbb{R}$ and
$\mathbb{H}^{n}\times \mathbb{R}$ in \cite{DCA}. In \cite{NC} Turgay
obtained complete classification of H-hypersurfaces with three
distinct principal curvatures in Euclidean spaces. Chen and Garray
in \cite{BY4} characterized $\delta(2)$-ideal null 2-type
hypersurfaces in Euclidean space as spherical cylinders. Also, Chen
has proved that every $\delta(3)$-ideal null $2$-type hypersurface
in Euclidean space has constant mean curvature \cite{BY5}.

For a Riemannian manifold $M^{n}$ with $n\geq 3$ and an integer $r
\in [2, n-1],$ Chen introduced the notion of $\delta$-invariant
$\delta(r)$ by

\begin{equation}\label{1.3}
\delta(r)(p)=\rho(p)-\inf_{r} \rho(L^{r}),
\end{equation}
where $\rho(p)$ is the scalar curvature at $p\in M^{n}$ and
$\rho(L^{r})$ is the scalar curvature of a linear subspace $L^{r}$
of dimension $r\geq 3$ of the tangent space $T_{p}(M)$.

For any n-dimensional submanifold $M^{n}$ in a Euclidean space
$\mathbb{E}^{m}$ and for an integer $r \in [2, n-1]$, Chen proved
the following universal sharp inequality \cite{BY6}

\begin{equation}\label{1.4}
\delta(r)(p)\leq \frac{n^{2}(n-r)}{2(n-r+1)}H^{2},
\end{equation}
where $H^{2}=\langle \vec{H},\vec{H}\rangle$ is the squared mean
curvature. If equality case for (\ref{1.4}) holds identically, then
$M^n$ is called $\delta(r)$-ideal submanifold in $\mathbb{E}^{m}$.

 Recently, the first author proved that
biconservative Lorentz hypersurfaces in $\mathbb{E}^{n+1}_{1}$
having complex eigenvalues must be of constant mean curvature
(\cite{DEEP}). In this paper, we study $\delta(2)$, $\delta(3)$ and
$\delta(4)$-ideal biconservative hypersurfaces in Euclidean space.
In Section $3$, we investigate $\delta(2)$-ideal biconservative
hypersurfaces in $\mathbb{E}^{n+1}$ as a generalization of
\cite[Theorem 3.2]{BY3} and we obtain the following result:

 \begin{thm}\label{Theorem3.1}
Every $\delta(2)$-ideal biconservative hypersurface in Euclidean
space $\mathbb{E}^{n+1}$ for $n\geq 3$ is minimal.
\end{thm}

 In Section $4$, we study $\delta(3)$-ideal biconservative
hypersurfaces in $\mathbb{E}^{5}$ and concluded the following
result:

\begin{thm}\label{Theorem4.1}
 Every $\delta(3)$-ideal biconservative hypersurface in Euclidean space
$\mathbb{E}^{5}$ has constant mean curvature.
\end{thm}

The above result can be considered as  generalization of the result
proved in \cite[Theorem 4.3]{BY3}. Finally, in Section $5$, we study
$\delta(4)$-ideal biconservative hypersurfaces in $\mathbb{E}^{6}$
with constant scalar curvature and obtain the following result:

  \begin{thm}\label{main}
Every $\delta(4)$-ideal biconservative hypersurface in Euclidean
space $\mathbb{E}^{6}$ with constant scalar curvature has constant
mean curvature.
\end{thm}

\section{\textbf{Preliminaries}}

   Let ($M^{n}, g$) be a hypersurface isometrically immersed in Euclidean space
$(\mathbb{E}^{n+1}, \overline g)$ and $g = \overline g_{|M}$.
  Let $\overline\nabla $ and $\nabla$ denote the linear connections on $\mathbb{E}^{n+1}$ and $M$, respectively. Then, the Gauss and Weingarten formulae are given by
\begin{equation}\label{2.1}
\overline\nabla_{X}Y = \nabla_{X}Y + h(X, Y), \hspace{.2 cm} \forall
\hspace{.2 cm}X, Y \in\Gamma(TM),
\end{equation}
\begin{equation}\label{2.2}
\overline\nabla_{X}\xi = -\mathcal{A}_{\xi}X,
\end{equation}
where $\xi$ be the unit normal vector to $M$, $h$ is the second
fundamental form and $\mathcal{A}$ is the shape operator. It is well
known that the second fundamental form $h$ and shape operator
$\mathcal{A}$ are related by
\begin{equation}\label{2.3}
\overline{g}(h(X,Y), \xi) = g(\mathcal{A}_{\xi}X,Y).
\end{equation}

 The mean curvature is given by
\begin{equation}\label{2.4}
H = \frac{1}{n} \mbox{trace} \mathcal{\mathcal{A}}.
\end{equation}

 The Gauss and Codazzi equations are given by
\begin{equation}\label{2.5}
R(X, Y)Z = g(\mathcal{A}Y, Z) \mathcal{A}X - g(\mathcal{A}X, Z)
\mathcal{A}Y,
\end{equation}
\begin{equation}\label{2.6}
(\nabla_{X}\mathcal{A})Y = (\nabla_{Y}\mathcal{A})X,
\end{equation}
respectively, where $R$ is the curvature tensor and
\begin{equation}\label{2.7}
(\nabla_{X}\mathcal{A})Y = \nabla_{X}\mathcal{A}Y-
\mathcal{A}(\nabla_{X}Y)
\end{equation}
for all $ X, Y, Z \in \Gamma(TM)$.

The scalar curvature $\rho$ of $M$ is given by
\begin{equation}\label{2.8}
\rho = \frac{1}{2}(n^{2}H^{2}-\mbox{Trace} \mathcal{A}^{2}),
\end{equation}
\\

We need the following result from \cite[Theorem 13.3, 13.7]{BY2} (cf. also  corresponding propositions in \cite{BY3} and \cite{BY5}).\\

\begin{thm}\label{Theorem2.1}  Let $M^{n}$ be the hypersurface in Euclidean
space $\mathbb{E}^{n+1}$. Then for an integer $r \in [2, n-1]$
\begin{equation}\label{2.9}
\delta(r)(p)\leq \frac{n^{2}(n-r)}{2(n-r+1)}H^{2} ,
\end{equation}
 and equality holds at a point $p$ if and only
if there is an orthonormal basis
$\{e_{1},e_{2},e_{3},\ldots,e_{n}\}$ at $p$ such that the shape
operator is given by
\begin{center}
$ \mathcal{A} = \left(
                            \begin{array}{cccc}
                               D_r & 0   \\
                              0 & u_{r}I_{n-r} \\
                             \end{array}
                          \right),$ \hspace{.5 cm}
 \end{center}
where $D_{r}=$ ${\rm
diag}(\lambda_{1},\lambda_{2},\dots,\lambda_{r})$ and
$u_{r}=\lambda_{1}+\lambda_{2}+\dots+\lambda_{r}$ for some functions
$\lambda_{1}, \lambda_{2},\ldots \lambda_{r}$ defined on $M^{n}$.
 \end{thm}

\section{\textbf{$\delta(2)$-ideal biconservative hypersurfaces in $\mathbb{E}^{n+1}$}}

In this section we study $\delta(2)$-ideal biconservative
hypersurfaces in $\mathbb{E}^{n+1} (n > 2)$. From Theorem
\ref{Theorem2.1}, the shape operator for a $\delta(2)$-ideal
hypersurface in $E^{n+1}$ with respect to orthonormal basis
$\{e_{1},e_{2},\dots,e_{n}\}$ takes the form
\begin{equation}\label{4.1}
 \mathcal{A}= \left(
                            \begin{array}{cccccccccc}
                              \lambda_{1} &0 &0  & \cdots &0\\
                              0 & \lambda_{2} &0   & \cdots &0\\
                              0&0 &  \lambda_{1}+\lambda_{2}   &\cdots &0 \\
                              \vdots &\vdots &\vdots   &\ddots &\vdots \\
                              0 &0 &0   &\cdots &\lambda_{1}+\lambda_{2}\\
                            \end{array}
                          \right),
\end{equation}
for some functions $\lambda_{1}$, $\lambda_{2}$ defined on $M^{n}$,
which can be expressed as
\begin{equation}\label{4.2}
\mathcal{A}(e_{i})=\lambda_{i}e_{i}, \quad i=1,2,\dots,n,
\end{equation}
where $\lambda_{i}=\lambda_{1}+\lambda_{2}$ for $i=3,4,\dots,n.$

 Let us assume that the mean curvature is not constant and grad$H\neq 0$.
This implies the existence of an open connected subset $U$ of $M$,
with grad$_{p}H\neq 0$ for all $p\in U$. From (\ref{1.2}) it is easy
to see that grad$H$ is an eigenvector of the shape operator
$\mathcal{A}$ with corresponding principal curvature
$-\frac{nH}{2}$.

 Without lose of generality we choose $e_{1}$ in the direction of
grad$H$, which gives $\lambda_{1}=-\frac{nH}{2}$. We express grad$H$
as
\begin{equation}\label{3.2}
\mbox{grad} H =\sum_{i=1}^{n} e_{i}(H)e_{i}.
\end{equation}
 As we have taken $e_{1}$ parallel to grad$H$, it is
\begin{equation}\label{3}
e_{1}(H)\neq 0, \ \ e_{i}(H)= 0, \hspace{1 cm} i= 2, \dots, n.
\end{equation}
We express
\begin{equation}\label{3.4}
\nabla_{e_{i}}e_{j}=\sum_{k=1}^{n}\omega_{ij}^{k}e_{k}, \hspace{2
cm} i, j = 1, 2, \dots , n.
\end{equation}
Using (\ref{3.4}) and the compatibility conditions
$(\nabla_{e_{k}}g)(e_{i}, e_{i})= 0$, $(\nabla_{e_{k}}g)(e_{i},
e_{j})= 0$, we obtain
\begin{equation}\label{3.5}
\omega_{ki}^{i}=0, \hspace{1 cm} \omega_{ki}^{j}+ \omega_{kj}^{i}
=0,
\end{equation}
for $i \neq j, $ and $i, j, k = 1, 2,\dots, n$.\\

\emph{We  consider the following cases:}\\

\textbf{Case A.} $\lambda_{2}\neq \lambda_{A},\quad
A=3,4,\dots,n.$\\

 Taking $X=e_{i}, Y=e_{j}$, $(i\neq j)$ in (\ref{2.7}) and using (\ref{4.2}),
(\ref{3.4}), we get
$$(\nabla_{e_{i}}\mathcal{A})e_{j}=e_{i}(\lambda_{j})e_{j}+\sum_{k=1}^{n}\omega_{ij}^{k}e_{k}
(\lambda_{j}-\lambda_{k}).$$

Putting the value of $(\nabla_{e_{i}}\mathcal{A})e_{j}$ in
(\ref{2.6}), we find
$$e_{i}(\lambda_{j})e_{j}+\sum_{k=1}^{n}\omega_{ij}^{k}e_{k}
(\lambda_{j}-\lambda_{k})=e_{j}(\lambda_{i})e_{i}+\sum_{k=1}^{n}\omega_{ji}^{k}e_{k}
(\lambda_{i}-\lambda_{k}),$$ whereby taking inner product with
$e_{j}$ and $e_{k}$, we obtain
\begin{equation}\label{3.6}
e_{i}(\lambda_{j})=
(\lambda_{i}-\lambda_{j})\omega_{ji}^{j}=(\lambda_{j}-\lambda_{i})\omega_{jj}^{i},
\end{equation}
\begin{equation}\label{3.7}
(\lambda_{j}-\lambda_{k})\omega_{ij}^{k}=
(\lambda_{i}-\lambda_{k})\omega_{ji}^{k},
\end{equation}
respectively, for distinct $i, j, k = 1, 2, \dots, n.$

 Using
(\ref{3}), (\ref{3.4}) and the fact that $[e_{i} \hspace{.1 cm}
e_{j}](H)=0=\nabla_{e_{i}}e_{j}(H)-\nabla_{e_{j}}e_{i}(H)=\omega_{ij}^{1}e_{1}(H)-\omega_{ji}^{1}e_{1}(H),$
for $i\neq j$ and $i, j=2, \dots, n$, we find
\begin{equation}\label{3.9}
\omega_{ij}^{1}=\omega_{ji}^{1}.
\end{equation}

Using (\ref{2.4}), (\ref{4.1}) and $\lambda_{1}=-\frac{nH}{2}$, we
obtain
\begin{equation}\label{4.3}
\lambda_{2}=\frac{n(n+1)}{2(n-1)}H, \quad \lambda_{A}=\frac{nH}{n-1}
\quad A=3,4,\dots, n.
\end{equation}

Therefore, using (\ref{3}) and (\ref{4.3}), we obtain
\begin{equation}\label{4.4}
e_{1}(\lambda_{i})\neq 0,\quad e_{j}(\lambda_{i})=0,
\end{equation}
for $i=1,2,\dots, n$ and $j=2,3,4,\dots, n.$

Now, it can be seen that $\lambda_{1}$ can never be equal to
$\lambda_{2}$ and $\lambda_{A}$ for $A=  3, 4,\dots, n.$
 Indeed, if $\lambda_{1}= \lambda_{2}$ or $\lambda_{A}$, from (\ref{3.6}), we
 find
\begin{equation}\label{3.10}
 e_{1}(\lambda_{j})= (\lambda_{1}-\lambda_{j})\omega_{j1}^{j}=0,
 \quad j=2,A,
\end{equation}
which contradicts the first expression of (\ref{4.4}).

 Putting
$i\neq 1, j = 1, 2, A$ in (\ref{3.6}) and using (\ref{4.4}) and
(\ref{3.5}), we
 find
\begin{equation}\label{4.5}
\omega_{1i}^{1}=\omega_{2i}^{2}=
\omega_{Ai}^{A}=\omega_{11}^{i}=\omega_{22}^{i}= \omega_{AA}^{i}= 0,
\hspace{1 cm} i= 1, 2, A.
\end{equation}

Putting $k = 1, $ and $ i=2$, $j = A$ in (\ref{3.7}), and using
(\ref{3.5}), we
 get
\begin{equation}\label{4.6}
\omega_{2A}^{1}=\omega_{A2}^{1} = \omega_{21}^{A}= \omega_{A1}^{2}=
0.
\end{equation}

Now, putting $k = A,$ and $ i=\widetilde{A}$, $j = 1,2$ in
(\ref{3.7}), and using (\ref{3.5}), we
 get
\begin{equation}\label{4.7}
\omega_{\widetilde{A}1}^{A}=\omega_{A\widetilde{A}}^{1} =
\omega_{\widetilde{A}2}^{A}=\omega_{A\widetilde{A}}^{2} = 0,
\end{equation}
where $A\neq \widetilde{A}$ and $A,\widetilde{A}=3,4\dots,n.$

Now, evaluating $g(R(X,Y)Z,W)$, using (\ref{4.5})$\sim$(\ref{4.7}),
Gauss equation
(2.5) and (\ref{4.3}), we obtain the following:\\

 $\bullet\hspace{.2 cm}
 \mbox{For}\hspace{.2 cm} X=e_{1}, Y=e_{i}, Z=e_{1}, W=e_{i},$
\begin{equation}\label{4.8}
 e_{1}(\omega_{ii}^{1})- (\omega_{ii}^{1})^{2}= -\frac{nH}{2}
 \lambda_{i}, \quad i=2,3,\dots,n.
\end{equation}

 $\bullet\hspace{.2 cm}
 \mbox{For}\hspace{.2 cm} X=e_{2}, Y=e_{A}, Z=e_{2}, W=e_{A},$
\begin{equation}\label{4.9}
 \omega_{22}^{1}\omega_{AA}^{1}= -\frac{n^{2}(n+1)}{2(n-1)^{2}}H^{2}, \quad A=3,4,\dots,n.
\end{equation}

Now, putting $i=1$ and $j=2,A$ in (\ref{3.6}) and using (\ref{4.3}),
we find
\begin{equation}\label{4.9}
 e_{1}(H)= \frac{2nH}{n+1}\omega_{22}^{1},
\end{equation}
\begin{equation}\label{4.10}
 e_{1}(H)= \frac{(n+1)H}{2}\omega_{AA}^{1},\quad A=3,4,\dots,n,
\end{equation}
respectively.

Equating (\ref{4.9}) and (\ref{4.10}), we get
\begin{equation}\label{4.11}
 \omega_{22}^{1}= \frac{(n+1)^{2}}{4n}\omega_{AA}^{1},\quad A=3,4,\dots,n,
\end{equation}
which by using (\ref{4.9}) gives
\begin{equation}\label{4.12}
 (\omega_{22}^{1})^{2}= -\frac{n(n+1)^{3}}{8(n-1)^{2}}H^{2},
\end{equation}
\begin{equation}\label{4.13}
 (\omega_{AA}^{1})^{2}= -\frac{2n^{3}}{(n+1)(n-1)^{2}}H^{2},\quad A=3,4,\dots,n,
\end{equation}

Now, differentiating (\ref{4.11}) along $e_{1}$ and using
(\ref{4.8}), (\ref{4.12}) and (\ref{4.13}), we obtain
\begin{equation}\label{4.14}
 H^{2}(n-1)^{2}=0,
\end{equation}
which gives $H=0$ as $n> 2.$\\

\textbf{Case B.} $\lambda_{2}= \lambda_{A},\quad
A=3,4,\dots,n.$\\

In this case, using (\ref{4.3}), we obtain that
$H[\frac{n(n+1)}{2(n-1)}-\frac{n}{n-1}]=0,$ which implies
$n(n-1)H=0.$
Since $n>2$, it is $H=0.$\\

 Combining cases \textbf{A} and \textbf{B}, we can obtain
 Theorem \ref{Theorem3.1}.\\

\section{\textbf{$\delta(3)$-ideal biconservative hypersurfaces in $\mathbb{E}^{5}$}}

In this section we study $\delta(3)$-ideal biconservative
hypersurfaces in $\mathbb{E}^{5}$. From Theorem \ref{Theorem2.1},
the shape operator for a $\delta(3)$-ideal hypersurface in
$\mathbb{E}^{5}$ with respect to orthonormal basis
$\{e_{1},e_{2},e_{3},e_{4}\}$ takes the form
\begin{equation}\label{5.1}
 \mathcal{A}= \left(
                            \begin{array}{cccccccc}
                              \lambda_{1} & & & \\
                              &   \lambda_{2} & & \\
                              & &     \lambda_{3} & \\
                              & & & \lambda_{1}+\lambda_{2}+\lambda_{3} \\
                            \end{array}
                          \right).
\end{equation}
for some functions $\lambda_{1}$, $\lambda_{2}$, $\lambda_{3}$
defined on $M^{4}$, which can be expressed as
\begin{equation}\label{5.2}
\mathcal{A}(e_{i})=\lambda_{i}e_{i}, \quad i=1,2,3,4,
\end{equation}
where $\lambda_{4}=\lambda_{1}+\lambda_{2}+\lambda_{3}.$

 Let us assume that the mean curvature is not constant and grad$H\neq 0$.
This implies the existence of an open connected subset $U$ of $M$
with grad$_{p}H\neq 0$, for all $p\in U$. From (\ref{1.2}) it is
easy to see that grad$H$ is an eigenvector of the shape operator
$\mathcal{A}$ with corresponding principal curvature $-2H$.

Now, if grad$H$ is in the direction of $e_{4}$ then
$\lambda_{1}+\lambda_{2}+\lambda_{3}=-2H$. Since from (\ref{2.4}) we
have $2(\lambda_{1}+\lambda_{2}+\lambda_{3})=4H$, this implies that
$H=0$, which is a contradiction. Therefore, without loss of
generality we may choose $e_{1}$ in the direction of grad$H$, which
gives $\lambda_{1}=-2H.$ As grad$H$ =
$e_{1}(H)e_{1}+e_{2}(H)e_{2}+e_{3}(H)e_{3}+e_{4}(H)e_{4}$,  we have
\begin{equation}\label{5.3}
e_{1}(H)\neq 0, e_{i\ \ }(H)= 0, \hspace{1 cm} i= 2, 3, 4.
\end{equation}

Also, in this section, equations (\ref{3.4}), (\ref{3.5}),
(\ref{3.6}), (\ref{3.7}) and (\ref{3.9}) hold for $n=4$.

Now, we can show that $\lambda_{j}\neq \lambda_{1}, j= 2, 3, 4,$ in
a similar way as we have shown in Section $3$.\\

\emph{We  consider the following cases:}\\

\textbf{Case A.} $\lambda_{2}\neq \lambda_{3}\neq
\lambda_{4}$\\

Using $\lambda_{1}=-2H$ and equation (\ref{2.4}), we obtain that
$\lambda_{4}=2H$ and
 \begin{equation}\label{5.11}
 \lambda_{2}+ \lambda_{3}= 4H.
\end{equation}

 Putting $i\neq 1, j = 1, 4$ in (\ref{3.6}) and using (\ref{5.3}) and (\ref{3.5}), we
 find
\begin{equation}\label{5.12}
\omega_{1i}^{1}=  \omega_{4j}^{4}=\omega_{11}^{i}=\omega_{44}^{j}=
0, \hspace{1 cm} j=2,3,4,\quad i= 1, 2, 3, 4.
\end{equation}

Putting $k = 1, j\neq i, $ and $ i, j = 2, 3, 4$ in (\ref{3.7}), and
using (\ref{3.9}), we
 get
\begin{equation}\label{5.13}
\omega_{ij}^{1}=\omega_{ji}^{1} = \omega_{1i}^{j}= \omega_{i1}^{j}=
0, \hspace{1 cm} j\neq i, \mbox{and} \hspace{.2 cm} i, j = 2, 3, 4.
\end{equation}

Thus, we have the following:

\begin{lem}\label{lemma3.11} Let $M^{4}$ be a $\delta(3)$-ideal
biconservative hypersurface of non constant mean curvature in
Euclidean space $\mathbb{E}^{5}$. Then, we obtain
\begin{equation}\label{5.14}
 \nabla_{e_{1}}e_{i}= 0, \hspace{.2 cm} i= 1, 2, 3, 4
\end{equation}
\begin{equation}\label{5.15}
 \nabla_{e_{4}}e_{4}=\omega_{44}^{1}e_{1}, \quad \nabla_{e_{i}}e_{1}=\omega_{i1}^{i}e_{i},  \hspace{.2 cm} i= 2, 3,
 4,
\end{equation}
\begin{equation}\label{5.16}
\nabla_{e_{i}}e_{4}=\sum_{ k=2}^{3}\omega_{i4}^{k}e_{k},  \quad
 \nabla_{e_{i}}e_{i}=
\sum_{i\neq j, j=1}^{4}\omega_{ii}^{j}e_{j}, \hspace{.2 cm} i= 2, 3,
\end{equation}
\begin{equation}\label{5.17}
 \nabla_{e_{4}}e_{j}= \sum_{k\neq j, k=2}^{3}\omega_{4j}^{k}e_{k},\quad \nabla_{e_{i}}e_{j}= \sum_{k\neq j, k=2}^{4}\omega_{ij}^{k}e_{k},\hspace{.2 cm} i,j= 2, 3, and \hspace{.2 cm}i\neq j,
\end{equation}
where $\omega_{jk}^{i}$ satisfy (\ref{3.5}),  (\ref{3.6}) and
(\ref{3.7}).
\end{lem}

Evaluating $g( R(X,Y)Z,W)$, using Lemma \ref {lemma3.11} and Gauss
equation
(\ref{2.5}), we find the following:\\

  $\bullet\hspace{.2 cm}
 \mbox{For}\hspace{.2 cm} X=e_{1}, Y=e_{i}, Z=e_{1}, W=e_{i},$
\begin{equation}\label{5.18}
 e_{1}(\omega_{ii}^{1})- (\omega_{ii}^{1})^{2}= -2H
 \lambda_{i}, \quad i=2,3,4.
\end{equation}

  $\bullet \hspace{.2 cm}
 \mbox{For}\hspace{.2 cm}
X=e_{1}, Y=e_{i}, Z=e_{i}, W=e_{j},$
\begin{equation}\label{5.19}
 e_{1}(\omega_{ii}^{j})- \omega_{ii}^{j} \omega_{ii}^{1}= 0, \quad
 i\neq j,\quad i,j=2,3,4.
\end{equation}

$\bullet \hspace{.2 cm}
 \mbox{For}\hspace{.2 cm} X=e_{1}, Y=e_{i}, Z=e_{j}, W= e_{4}$
 \begin{equation}\label{5.20}
 e_{1}(\omega_{ij}^{4})- \omega_{ii}^{1} \omega_{ij}^{4}= 0\quad
 i\neq j,\quad i,j=2,3.
\end{equation}
\\

Now, we have:
%
\begin{lem}\label{lemma3.22} Let $M^{4}$ be a $\delta(3)$-ideal
biconservative hypersurface of non constant mean curvature in
Euclidean space $\mathbb{E}^{5}$. Then,
\begin{equation}\label{5.21}
\omega_{ij}^{4}= 0,
\end{equation}
 for $i,j $ = 2, 3.
\end{lem}
\noindent {\it Proof.} Differentiating (\ref{5.11}) along $e_{4}$
and using (\ref{3.6}) and $\lambda_{4}=2H$, we obtain
\begin{equation}\label{5.22}
(\lambda_{2}-2H)\omega_{22}^{4}+(\lambda_{3}-2H)\omega_{33}^{4}= 0.
\end{equation}

Now, differentiating (\ref{5.22}) along $e_{1}$ and using
(\ref{5.19}) and $\lambda_{1}=-2H$,, we get
\begin{equation}\label{5.23}
(\lambda_{2}\omega_{22}^{1}-e_{1}(H))\omega_{22}^{4}+(\lambda_{3}\omega_{33}^{1}-e_{1}(H))\omega_{33}^{4}=0.
\end{equation}

Equations (\ref{5.22}) and (\ref{5.23}) consitute  a homogeneous
system in two variables $\omega_{22}^{4}$ and $\omega_{33}^{4}$,
having either non trivial solution or trivial solution. If it has
trivial solution only,  then  $\omega_{22}^{4}=0$ and
$\omega_{33}^{4}=0$.

We will show that it is not possible to have a non trivial solution.
Indeed, if it had one,  then the determinant formed by the
coefficients of $\omega_{22}^{4}$ and $\omega_{33}^{4}$ in
(\ref{5.22}) and (\ref{5.23}) would be zero, i.e.
\begin{equation}\label{5.24}
(\lambda_{3}-2H)(\lambda_{2}\omega_{22}^{1}-e_{1}(H))-(\lambda_{2}-2H)(\lambda_{3}\omega_{33}^{1}-e_{1}(H))=0.
\end{equation}

Differentiating (\ref{5.11}) along $e_{1}$ and using (\ref{3.6}) and
$\lambda_{1}=-2H$, we obtain
\begin{equation}\label{5.25}
(\lambda_{2}+2H)\omega_{22}^{1}+(\lambda_{3}+2H)\omega_{33}^{1}=
4e_{1}(H).
\end{equation}

Eliminating $e_{1}(H)$ from (\ref{5.24}) and (\ref{5.25}), we obtain
\begin{equation}\label{5.26}
(\lambda_{2}-2H)^{2}(\omega_{22}^{1}-\omega_{33}^{1})= 0.
\end{equation}
By our assumption it is $\lambda_{2}\neq 2H$.  If
$\omega_{22}^{1}=\omega_{33}^{1}$, then from (\ref{5.18}), we get
$\lambda_{2}=\lambda_{3}$, which is a contradiction to our
assumption. Hence, we have $\omega_{ii}^{4}=0$ for $i=2,3$.\\

Now, we want to prove that $\omega_{ij}^{4}=0$ for $i \neq j$,
$i,j=2,3$.

From (\ref{3.7}) and using $\lambda_{4}=2H$ we have
\begin{equation}\label{5.27}
(\lambda_{3}-2H)\omega_{23}^{4}=(\lambda_{2}-2H)\omega_{32}^{4}.
\end{equation}

Differentiating (\ref{5.27}) along $e_{1}$ and using (\ref{5.20})
and $\lambda_{1}=-2H$, we get
\begin{equation}\label{5.28}
[(\lambda_{3}+2H)\omega_{33}^{1}+(\lambda_{3}-2H)\omega_{22}^{1}-2e_{1}(H)]\omega_{23}^{4}=
[(\lambda_{2}+2H)\omega_{22}^{1}+(\lambda_{2}-2H)\omega_{33}^{1}-2e_{1}(H)]\omega_{32}^{4}.
\end{equation}

Now, equations (\ref{5.27}) and (\ref{5.28}) constitute a
homogeneous system of equations in two variables $\omega_{23}^{4}$
and $\omega_{32}^{4}$, having either non trivial solution or trivial
solution only. If it has trivial solution only, then
$\omega_{23}^{4}=0$ and $\omega_{32}^{4}=0$.

If it has non trivial solution also, then the determinant formed by
the coefficients of $\omega_{23}^{4}$ and $\omega_{32}^{4}$ in
(\ref{5.27}) and (\ref{5.28}) will be zero, i.e.
\begin{equation}\label{5.29}
2H(\lambda_{3}-2H)\omega_{22}^{1}-2H(\lambda_{2}-2H)\omega_{33}^{1}+e_{1}(H)(\lambda_{2}-\lambda_{3})=0.
\end{equation}

Eliminating $e_{1}(H)$ from (\ref{5.29}), using (\ref{5.25}), we
obtain
\begin{equation}\label{5.30}
(\lambda_{2}-2H)^{2}(\omega_{22}^{1}-\omega_{33}^{1})= 0,
\end{equation}
which is not possible from (\ref{5.26}). Hence, $\omega_{ij}^{4}=0$
for $i \neq j$, $i,j=2,3$, which proves  Lemma \ref{lemma3.22}. \qed

\smallskip
We now evaluate $g(R(X,Y)Z,W)$ using Lemma \ref{lemma3.11}, Lemma
\ref{lemma3.22}, Gauss equation
(\ref{2.5}) and  $\lambda_{4}=2H$.  We obtain the following:\\

  $\bullet\hspace{.2 cm}
 \mbox{For}\hspace{.2 cm} X=e_{2}, Y=e_{4}, Z=e_{2}, W= e_{4}$
\begin{equation}\label{5.31}
\begin{array}{rcl}
-\omega_{22}^{1}\omega_{44}^{1} = \lambda_{2} \lambda_{4}.
\end{array}
\end{equation}

 $\bullet\hspace{.2 cm}
 \mbox{For}\hspace{.2 cm} X=e_{3}, Y=e_{4}, Z=e_{3}, W= e_{4}$
\begin{equation}\label{5.32}
\begin{array}{rcl}
-\omega_{33}^{1}\omega_{44}^{1} = \lambda_{3} \lambda_{4}.
\end{array}
\end{equation}

Now, using $\lambda_{1}=-2H$, $\lambda_{4}=2H$, and (\ref{3.6}), we
get
\begin{equation}\label{5.33}
e_{1}(H)=2H \omega^{1}_{44},
\end{equation}
which by differentiating again along $e_{1}$ gives
\begin{equation}\label{5.34}
e_{1}e_{1}(H)=\frac{3e_{1}^{2}(H)}{2H}-8H^{3},
\end{equation}

Adding (\ref{5.31}) and (\ref{5.32}) and using (\ref{5.33}) and
(\ref{5.11}), we obtain

\begin{equation}\label{5.35}
(\omega_{22}^{1}+\omega_{33}^{1})e_{1}(H) = -16 H^{3}.
\end{equation}

By solving (\ref{5.25}) and (\ref{5.35}) for $\omega^{1}_{22}$ and
$\omega^{1}_{33}$ and using (\ref{5.11}), we find
\begin{equation}\label{5.36}\begin{array}{rcl}
\omega_{22}^{1}=
-\Big[\frac{8H^{3}(2H+\lambda_{3})}{e_{1}(H)(2H-\lambda_{3})}+\frac{2e_{1}^{2}(H)}{e_{1}(H)(2H-\lambda_{3})}\Big].\end{array}
\end{equation}
\begin{equation}\label{5.37}\begin{array}{rcl}
\omega_{33}^{1}=
-\Big[\frac{8H^{3}(6H-\lambda_{3})}{e_{1}(H)(2H-\lambda_{3})}+\frac{2e_{1}^{2}(H)}{e_{1}(H)(2H-\lambda_{3})}\Big].\end{array}
\end{equation}

By squaring and adding (\ref{5.36}), (\ref{5.37}), we obtain
\begin{equation}\label{5.38}\begin{array}{rcl}
(\omega_{22}^{1})^{2}+(\omega_{33}^{1})^{2}=
\frac{8}{e_{1}^{2}(H)(2H-\lambda_{3})^{2}}[16H^{6}(\lambda_{3}^{2}+20H^{2}-4H\lambda_{3})+e_{1}^{4}(H)\\+32H^{4}e_{1}^{2}(H)].
\end{array}
\end{equation}

Differentiating (\ref{5.35}) along $e_{1}$, and using (\ref{5.18})
and (\ref{5.11}), we get
\begin{center}
$[-8H^{2}+(\omega_{22}^{1})^{2}+(\omega_{33}^{1})^{2}]e_{1}(H)+e_{1}e_{1}(H)(\omega_{22}^{1}+\omega_{33}^{1})=-48H^{3},$
\end{center}
which by using (\ref{5.34}), (\ref{5.35}) and (\ref{5.38}) gives
\begin{equation}\label{5.39}\begin{array}{rcl}
f(e_{1}(H),\lambda_{3},H)=e_{1}^{4}(H)+2H^{2}(\lambda_{3}^{2}+20H^{2}-4H\lambda_{3})e_{1}^{2}(H)\\+36H^{6}(\lambda_{3}^{2}+12H^{2}-4H\lambda_{3})=0.
\end{array}
\end{equation}

 If we differentiate (\ref{5.39}) along $e_{1}$ and using
(\ref{3.6}), (\ref{5.34}) and (\ref{5.37}), we obtain a polynomial
\begin{equation}\label{5.40}\begin{array}{rcl}
g(e_{1}(H),\lambda_{3},H)=0.
\end{array}
\end{equation}

 We rewrite $f(e_{1}(H),\lambda_{3},H)$, $g(e_{1}(H),\lambda_{3},H)$ as polynomials
$f_{(H, \lambda_{3})}(e_{1}(H))$ and $g_{(H,
\lambda_{3})}(e_{1}(H))$ of $e_{1}(H)$ with coefficients in
polynomial ring $R_{1}[H,\lambda_{3}]$ over real field $\mathbb{R}$.
We know that equations $f_{(H, \lambda_{3})}(e_{1}(H)) = 0$ and
$g_{(H, \lambda_{3})}(e_{1}(H)) = 0$ have a common root if and only
if resultant $\Re(f_{(H, \lambda_{3})}, g_{(H, \lambda_{3})}) = 0$.
It is obvious that $\Re(f_{(H,\lambda_{3})}, g_{(H, \lambda_{3})})$
is a polynomial of $\lambda_{3}$ and $H$. So, we have
\begin{equation}\label{5.41}
\widetilde{f}(\lambda_{3}, H)=\Re(f_{(H, \lambda_{3})}, g_{(H,
\lambda_{3})}) = 0
\end{equation}

If we differentiate (\ref{5.41}) along $e_{1}$ and using (\ref{3.6})
and (\ref{5.37}), we obtain again a polynomial
\begin{equation}\label{5.42}
\widetilde{g}(\lambda_{3},H)=0.
\end{equation}

Again, we rewrite $\widetilde{f}(\lambda_{3},H)$,
$\widetilde{g}(\lambda_{3},H)$ as polynomials
$\widetilde{f}_{H}(\lambda_{3}), \widetilde{g}_{H}(\lambda_{3})$ of
$\lambda_{3}$ with coefficients in the polynomial ring $R[H]$ over
$\mathbb{R}$. Since $\widetilde{f}_{H}(\lambda_{3}) =
\widetilde{g}_{H}(\lambda_{3}) = 0$ and $\lambda _3$ is a common
root of $\widetilde{f}_H, \widetilde{g}_H$, hence resultant
$\Re(\widetilde{f}_{H}, \widetilde{g}_{H}) = 0$. It is obvious that
$\Re(\widetilde{f}_{H}, \widetilde{g}_{H})$ is a polynomial
of $H$ with constant coefficients, therefore $H$ must be a constant.\\

\textbf{Case B.}  $\lambda_{2}= \lambda_{3}.$\\

In this case, using (\ref{5.11}), we find $\lambda_{2}=
\lambda_{3}=2H=\lambda_{4}$, which by using (\ref{3.6}) gives
$\omega_{22}^{1}=\omega_{33}^{1}=\omega_{44}^{1}.$

Now, it can be easily seen that our equation (\ref{5.31}) reduces to
$(\omega_{22}^{1})^{2}=-4H^{2}$, which is possible only if $H=0$,
since
$\omega_{22}^{1}$ is a real connection coefficient.\\

Hence, combining cases \textbf{A} and \textbf{B}, we can conclude
Theorem \ref{Theorem4.1}.\\

\section{\textbf{$\delta(4)$-ideal biconservative hypersurfaces in $\mathbb{E}^{6}$}}

In this section we  study $\delta(4)$-ideal biconservative
hypersurfaces in $\mathbb{E}^{6}$ having constant scalar curvature.
From Theorem \ref{Theorem2.1} the shape operator for a
$\delta(4)$-ideal hypersurface in $\mathbb{E}^{6}$ with respect to
orthonormal basis $\{e_{1},e_{2},e_{3},e_{4},e_{5}\}$ takes the form
\begin{equation}\label{3.1}
 \mathcal{A}= \left(
                            \begin{array}{cccccccccc}
                              \lambda_{1} & & & &\\
                              &   \lambda_{2} & & &\\
                              & &     \lambda_{3} & &\\
                              & & & \lambda_{4} &\\
                              & & & & \lambda_{1}+\lambda_{2}+\lambda_{3}+\lambda_{4} \\
                            \end{array}
                          \right),
\end{equation}
for some functions $\lambda_{1}$, $\lambda_{2}$, $\lambda_{3}$,
$\lambda_{4}$ defined on $M^{5}$, which can be expressed as
\begin{equation}\label{3.2}
\mathcal{A}(e_{i})=\lambda_{i}e_{i}, \quad i=1,2,3,4,5,
\end{equation}
where $\lambda_{5}=\lambda_{1}+\lambda_{2}+\lambda_{3}+\lambda_{4}.$

 We assume that the mean curvature is not constant and grad$H\neq 0$.
This implies the existence of a open connected subset $U$ of $M$,
with grad$_{p}H\neq 0$ for all $p\in U$. From (\ref{1.2}), it is
easy to see that grad$H$ is an eigenvector of the shape operator
$\mathcal{A}$ with the corresponding principal curvature
$-\frac{5H}{2}$.

If grad$H$ is in the direction of $e_{5}$ then
$\lambda_{1}+\lambda_{2}+\lambda_{3}+\lambda_{4}=-\frac{5H}{2}$.
Since from (\ref{2.4}), we have
$2(\lambda_{1}+\lambda_{2}+\lambda_{3}+\lambda_{4})=5H$ which
implies $H=0$. This is a contradiction. Therefore, without losing
generality, we choose $e_{1}$ in the direction of grad$H$, which
gives $\lambda_{1}=-\frac{5H}{2}$ and we have
\begin{equation}\label{3.3}
e_{1}(H)\neq 0, e_{i}(H)= 0, \hspace{2 cm} i= 2, 3, 4, 5.
\end{equation}

In this section also, equations (\ref{3.4}), (\ref{3.5}),
(\ref{3.6}), (\ref{3.7}) and (\ref{3.9}) hold true for $n=5$.

We can show that $\lambda_{j}\neq \lambda_{1}, j= 2, 3, 4, 5$ in
similar way as we have shown in Section $3$.\\

\emph{We consider the following cases:}\\

\textbf{Case A.} $\lambda_{2}, \lambda_{3}, \lambda_{4}$ and
$\lambda_{5}$ are all distinct.\\

 Using $\lambda_{1}=-\frac{5H}{2}$ and (\ref{3.3}), we obtain that
 $\lambda_{5}=\frac{5H}{2}$ and
\begin{equation}\label{3.8}
e_{1}(\lambda_{1})\neq 0,\quad e_{i}(\lambda_{1})= 0,\quad
e_{i}(\lambda_{5}) = 0, \quad i = 2, 3, 4, 5.
\end{equation}

Using $\lambda_{1}=-\frac{5H}{2}$, $\lambda_{5}=\frac{5H}{2}$ and
equation (\ref{2.4}), we obtain that
 \begin{equation}\label{3.11}
\lambda_{2}+ \lambda_{3}+\lambda_{4}= 5H.
\end{equation}

 Putting $i\neq 1, j = 1, 5$ in (\ref{3.6}) and using (\ref{3.8}) and (\ref{3.5}), we
 find
\begin{equation}\label{3.12}
\omega_{1i}^{1}= \omega_{5j}^{5}=\omega_{11}^{i}=\omega_{55}^{j}= 0,
\hspace{1 cm} j=2,3,4,5,\quad i= 1, 2, 3, 4, 5.
\end{equation}

Putting $k = 1, j\neq i, $ and $ i, j = 2, 3, 4, 5$ in (\ref{3.7}),
and using (\ref{3.9}), we
 get
\begin{equation}\label{3.13}
\omega_{ij}^{1}=\omega_{ji}^{1} = \omega_{1i}^{j}= \omega_{i1}^{j}=
0.
\end{equation}

Thus, using (\ref{3.12}) and \ref{3.13}), we have the following:

\begin{lem}\label{lemma3.1} Let $M^{5}$ be a $\delta(4)$-ideal
biconservative hypersurface of non constant mean curvature in
Euclidean space $\mathbb{E}^{6}$. Then, we obtain
\begin{equation}\label{3.14}
 \nabla_{e_{1}}e_{i}= 0, \hspace{.2 cm} i= 1, 2, 3, 4, 5
\end{equation}
\begin{equation}\label{3.15}
 \nabla_{e_{5}}e_{5}=\omega_{55}^{1}e_{1}, \quad \nabla_{e_{i}}e_{1}=\omega_{i1}^{i}e_{i},  \hspace{.2 cm} i= 2, 3,
 4, 5,
\end{equation}
\begin{equation}\label{3.16}
\nabla_{e_{i}}e_{5}=\sum_{ k=2}^{4}\omega_{i5}^{k}e_{k},  \quad
 \nabla_{e_{i}}e_{i}=
\sum_{i\neq j, j=1}^{5}\omega_{ii}^{j}e_{j}, \hspace{.2 cm} i= 2, 3,
4,
\end{equation}
\begin{equation}\label{3.17}
 \nabla_{e_{5}}e_{j}= \sum_{k\neq j, k=2}^{4}\omega_{5j}^{k}e_{k},\quad \nabla_{e_{i}}e_{j}= \sum_{k\neq j, k=2}^{5}\omega_{ij}^{k}e_{k},\hspace{.2 cm} i,j= 2, 3, 4, and \hspace{.2 cm}i\neq j,
\end{equation}
where $\omega_{jk}^{i}$ satisfy (\ref{3.5}),  (\ref{3.6}) and
(\ref{3.7}).
\end{lem}

Evaluating $g( R(X,Y)Z,W)$, using Lemma \ref {lemma3.1} and Gauss
equation
(\ref{2.5}), we find the following:\\

  $\bullet\hspace{.2 cm}
 \mbox{For}\hspace{.2 cm} X=e_{1}, Y=e_{i}, Z=e_{1}, W=e_{i},$
\begin{equation}\label{3.18}
 e_{1}(\omega_{ii}^{1})- (\omega_{ii}^{1})^{2}= -\frac{5H}{2}
 \lambda_{i}, \quad i=2,3,4,5.
\end{equation}

  $\bullet \hspace{.2 cm}
 \mbox{For}\hspace{.2 cm}
X=e_{1}, Y=e_{i}, Z=e_{i}, W=e_{j},$
\begin{equation}\label{3.19}
 e_{1}(\omega_{ii}^{j})- \omega_{ii}^{j} \omega_{ii}^{1}= 0, \quad
 i\neq j,\quad i,j=2,3,4,5.
\end{equation}

$\bullet \hspace{.2 cm}
 \mbox{For}\hspace{.2 cm} X=e_{i}, Y=e_{j}, Z=e_{i}, W= e_{1}$
 \begin{equation}\label{3.20}
 e_{j}(\omega_{ii}^{1})+ \omega_{ii}^{j} \omega_{jj}^{1}-\omega_{ii}^{j} \omega_{ii}^{1}= 0\quad
 i\neq j,\quad i,j=2,3,4,5.
\end{equation}
\\

Now, using $\lambda_{1}=-\frac{5H}{2}$, $\lambda_{5}=\frac{5H}{2}$,
and (\ref{3.6}), we get
\begin{equation}\label{3.21}
e_{1}(H)=2H \omega^{1}_{55}.
\end{equation}

Differentiating (\ref{3.21}) along $e_{1}$ and using (\ref{3.18})
for $i=5$, we get
\begin{equation}\label{3.51}
2He_{1}e_{1}(H)= 3e_{1}^{2}(H)-25H^{4}.
\end{equation}

From (\ref{3.1}) and Gauss equation, the scalar curvature $\rho$
(constant) is given by
\begin{equation}\label{3.22}
 \rho= \frac{25H^{2}}{2}- \lambda_{2}^{2}- \lambda_{3}^{2}- \lambda_{4}^{2}.
\end{equation}

Now, we have:
%
\begin{lem}\label{lemma3.3} Let $M^{5}$ be a $\delta(4)$-ideal
biconservative hypersurface of non constant mean curvature in
Euclidean space $\mathbb{E}^{6}$. Then,
\begin{equation}\label{3.36}
\omega_{ij}^{5}= 0,
\end{equation}
 for $i,j $ = 2, 3, 4.
\end{lem}
\noindent {\it Proof.} Using (\ref{3.22}) and (\ref{3.11}), we get
\begin{equation}\label{3.37}
\lambda_{3}^{2}+\lambda_{4}^{2}+\lambda_{3}\lambda_{4}-5H(\lambda_{3}+\lambda_{4})=-\frac{25H^{2}}{4}-\frac{\rho}{2}.
\end{equation}

Differentiate (\ref{3.37}) along $e_{5}$ and using (\ref{3.6}), we
get
\begin{equation}\label{3.38}
\omega_{33}^{5}(2\lambda_{3}+\lambda_{4}-5H)(\lambda_{3}-\lambda_{5})+\omega_{44}^{5}(2\lambda_{4}+\lambda_{3}-5H)(\lambda_{4}-\lambda_{5})=
0,
\end{equation}
Again differentiating (\ref{3.38}) along $e_{1}$ and using
(\ref{3.19}) for $j=5$, $i=3,4$, we obtain
\begin{equation}\label{3.39}
\begin{array}{rcl}
\omega_{33}^{5}[2\omega_{33}^{1}\{(2\lambda_{3}+\lambda_{4}-5H)\lambda_{3}+(\lambda_{3}+\frac{5H}{2})(\lambda_{3}-\frac{5H}{2})\}
+\omega_{44}^{1}(\lambda_{4}+\frac{5H}{2})(\lambda_{3}-\\
\frac{5H}{2})-\frac{5}{2}e_{1}(H)(4\lambda_{3}+\lambda_{4}-10H)]
+\omega_{44}^{5}[2\omega_{44}^{1}\{(2\lambda_{4}+\lambda_{3}-5H)\lambda_{4}+(\lambda_{4}+\\
\frac{5H}{2})(\lambda_{4}-\frac{5H}{2})\}
+\omega_{33}^{1}(\lambda_{3}+\frac{5H}{2})(\lambda_{4}-\frac{5H}{2})-\frac{5}{2}e_{1}(H)(4\lambda_{4}+\lambda_{3}-10H)]=
0,\end{array}
\end{equation}

Now, we claim that $\omega^{5}_{33}=0$ and $\omega^{5}_{44}=0$.

Indeed, if $\omega^{5}_{33}\neq0$ and $\omega^{5}_{44}\neq0$, then
from (\ref{3.38}) and (\ref{3.39}), we have which gives
\begin{equation}\label{3.41}
f_{1}(\lambda_{3},\lambda_{4},H)\omega^{1}_{33}+g_{1}(\lambda_{3},\lambda_{4},H)\omega^{1}_{44}=h_{1}(\lambda_{3},\lambda_{4},H)e_{1}(H),
\end{equation}

where \begin{align*}\begin{array}{rcl}
f_{1}(\lambda_{3},\lambda_{4},H)&=2(3\lambda_{3}^{2}+\lambda_{3}\lambda_{4}-5H\lambda_{3}-\frac{25H^{2}}{4})
(4\lambda_{4}^{2}+2\lambda_{3}\lambda_{4}-5H\lambda_{3}-20H\lambda_{4}-25H^{2}),&\\
g_{1}(\lambda_{3},\lambda_{4},H)&=2(3\lambda_{4}^{2}+\lambda_{3}\lambda_{4}-5H\lambda_{4}-\frac{25H^{2}}{4})
(4\lambda_{3}^{2}+2\lambda_{3}\lambda_{4}-5H\lambda_{4}-20H\lambda_{3}-25H^{2}),&\\
h_{1}(\lambda_{3},\lambda_{4},H)&=5(\lambda_{3}-\lambda_{4})(-7\lambda_{3}\lambda_{4}+20H\lambda_{3}+20H\lambda_{4}-\frac{75}{2}H^{2}-2\lambda_{4}^{2}
-2\lambda_{3}^{2}-2\lambda_{3}\lambda_{4})
\end{array}
\end{align*}
are homogeneous polynomials in $\lambda_{3}$, $\lambda_{4}$ and $H$.

Differentiate (\ref{3.37}) along $e_{1}$ and using (\ref{3.6}), we
get
\begin{equation}\label{3.40}\begin{array}{rcl}
\omega_{33}^{1}(2\lambda_{3}+\lambda_{4}-5H)(\lambda_{3}-\lambda_{1})+\omega_{44}^{1}(2\lambda_{4}+\lambda_{3}-5H)(\lambda_{4}-\lambda_{1})\\=
[-\frac{25H}{2}+5(\lambda_{3}+\lambda_{4})]e_{1}(H),\end{array}
\end{equation}

Now, solving equation (\ref{3.41}) and (\ref{3.40}), it can be
easily seen that
\begin{equation}\label{3.42}
\omega^{1}_{33}=
\frac{f_{2}(\lambda_{3},\lambda_{4},H)}{\widetilde{f}_{2}(\lambda_{3},\lambda_{4},H)}e_{1}(H),
\quad \omega^{1}_{44}=
\frac{g_{2}(\lambda_{3},\lambda_{4},H)}{\widetilde{g}_{2}(\lambda_{3},\lambda_{4},H)}e_{1}(H)
\end{equation}
where $f_{2}(\lambda_{3},\lambda_{4},H)$,
$\widetilde{f}_{2}(\lambda_{3},\lambda_{4},H)$,
$g_{2}(\lambda_{3},\lambda_{4},H)$ and
$\widetilde{g}_{2}(\lambda_{3},\lambda_{4},H)$ are some homogeneous
polynomials in $\lambda_{3}$, $\lambda_{4}$ and $H$.

Differentiating (\ref{3.40}) along $e_{1}$ and using (\ref{3.6}),
(\ref{3.18}) and (\ref{3.51}), we get
\begin{equation}\label{3.43}\begin{array}{rcl}
(\omega_{33}^{1})^{2}(6\lambda_{3}^{2}+2\lambda_{3}\lambda_{4}-5H\lambda_{4}-\frac{25H^{2}}{2})+
(\omega_{44}^{1})^{2}(6\lambda_{4}^{2}+2\lambda_{3}\lambda_{4}-5H\lambda_{3}-\\
\frac{25H^{2}}{2})
+2\omega_{33}^{1}\omega_{44}^{1}(\lambda_{3}+\frac{5H}{2})(\lambda_{4}+\frac{5H}{2})
+\frac{5}{2}\omega_{33}^{1}e_{1}(H)(\lambda_{4}-2\lambda_{3}-15H)\\+\frac{5}{2}\omega_{44}^{1}e_{1}(H)(\lambda_{3}-2\lambda_{4}-15H)-
\frac{5}{2H}e_{1}^{2}(H)(3\lambda_{3}+3\lambda_{4}-20H)=
\frac{5H}{4}(4\lambda_{3}^{3}\\+4\lambda_{3}^{3}+10H\lambda_{3}\lambda_{4}-75H^{2}\lambda_{3}-75H^{2}\lambda_{4}
+2\lambda_{3}^{2}\lambda_{4}+2\lambda_{4}^{2}\lambda_{3}+125H^{3}),\end{array}
\end{equation}
which by eliminating $\omega_{33}^{1}$ and $\omega_{44}^{1}$ using
(\ref{3.42}) gives
\begin{equation}\label{3.44}
e_{1}^{2}(H)h_{2}(\lambda_{3},\lambda_{4},H)=\widetilde{h}_{2}(\lambda_{3},\lambda_{4},H).
\end{equation}

Now, differentiating (\ref{3.44}) along $e_{1}$ and using
(\ref{3.42}), (\ref{3.51}), we obtain
\begin{equation}\label{3.45}
e_{1}^{2}(H)h_{3}(\lambda_{3},\lambda_{4},H)=\widetilde{h}_{3}(\lambda_{3},\lambda_{4},H),
\end{equation}
where $h_{2}(\lambda_{3},\lambda_{4},H)$,
$\widetilde{h}_{2}(\lambda_{3},\lambda_{4},H)$,
$h_{3}(\lambda_{3},\lambda_{4},H)$ and
$\widetilde{h}_{3}(\lambda_{3},\lambda_{4},H)$ are homogeneous
polynomials in $\lambda_{3}$, $\lambda_{4}$ and $H$.

It can be easily seen that by eliminating $e_{1}^{2}(H)$ from
(\ref{3.44}) and (\ref{3.45}), we obtain a homogeneous polynomial
equation defined as
\begin{equation}\label{3.45a}
\alpha(\lambda_{3},\lambda_{4},H)=0,
\end{equation}
which by differentiating along $e_{1}$ and using (\ref{3.42}) gives
again a homogeneous polynomial equation defined as
\begin{equation}\label{3.45b}
\beta(\lambda_{3},\lambda_{4},H)=0.
\end{equation}

Again if we differentiate (\ref{3.45b}) along $e_{1}$ and using
(\ref{3.42}), we find a homogeneous polynomial equation
\begin{equation}\label{3.45c}
\gamma(\lambda_{3},\lambda_{4},H)=0.
\end{equation}

 We rewrite $\alpha(\lambda_{3},\lambda_{4},H)$, $\beta(\lambda_{3},\lambda_{4},H)$ as polynomials
$\alpha_{(H,\lambda_{4})}(\lambda_{3}),
\beta_{(H,\lambda_{4})}(\lambda_{3})$ of $\lambda_{3}$ with
coefficients in polynomial ring $R_{2}[\lambda_{4},H]$ over real
field $\mathbb{R}$. Also, Equations
$\alpha_{(H,\lambda_{4})}(\lambda_{3}) = 0$ and
$\beta_{(H,\lambda_{4})}(\lambda_{3}) = 0$ have a common root if and
only if resultant $\Re(\alpha_{(H,\lambda_{4})},
\beta_{(H,\lambda_{4})}) = 0$, which is a polynomial equation of
$\lambda_{4}$ and $H$ and can be defined as
\begin{equation}\label{3.45d}
f_{3}(\lambda_{4},H)=0.
\end{equation}

Similarly, we can eliminate $\lambda_{3}$ from (\ref{3.45b}) and
(\ref{3.45c}), we obtain another polynomial equation defined as
\begin{equation}\label{3.45e}
g_{3}(\lambda_{4},H)=0.
\end{equation}

Again, we rewrite $f_{3}(\lambda_{4},H)$, $g_{3}(\lambda_{4},H)$ as
polynomials $f_{3(H)}(\lambda_{4}), g_{3(H)}(\lambda_{4})$ of
$\lambda_{4}$ with coefficients in polynomial ring $R[H]$ over real
field $\mathbb{R}$. Since $f_{3(H)}(\lambda_{4}) = 0$ and
$g_{3(H)}(\lambda_{4}) = 0$ have a common root $\lambda_{4}$, which
gives resultant $\Re(f_{3(H)}, g_{3(H)}) = 0$. Clearly
$\Re(f_{3(H)}, g_{3(H)})$ is a polynomial of $H$ with constant
coefficients. So, $\Re(f_{3(H)}, g_{3(H)}) = 0$ which implies that
$H$ must be a constant which is contradiction to our assumption.
Hence, we have $\omega^{5}_{ii}= 0$ for $i=2,3,4$.
\\

Now, we claim that $\omega^{5}_{ij}= 0$ for $i\neq j$, $i,j=2,3,4.$

 Using Lemma \ref{lemma3.1} and (\ref{2.5}) to evaluate
$g(R(e_{1},e_{2})e_{3},e_{5})$, and $g(R(e_{1},e_{3})e_{2},e_{5})$,
we obtain
\begin{equation}\label{3.46}
e_{1}(\omega_{23}^{5})-\omega_{22}^{1}\omega_{23}^{5}=0,\quad\mbox{and}\quad
e_{1}(\omega_{32}^{5})-\omega_{33}^{1}\omega_{32}^{5}=0,
\end{equation}
respectively.

 Putting $j=5, k=2, i=3$ in (\ref{3.7}), we get
\begin{equation}\label{3.47}
(\lambda_{3}-\lambda_{5})\omega_{23}^{5} =
(\lambda_{2}-\lambda_{5})\omega_{32}^{5}.
\end{equation}

Differentiating (\ref{3.47}) with respect to $e_{1}$ and using
(\ref{3.6}), (\ref{3.46}), we get
\begin{equation}\label{3.48}
\begin{array}{rcl}
\omega_{32}^{5}[\omega_{33}^{1}(\lambda_{2}-\lambda_{5})+\omega_{22}^{1}(\lambda_{2}-\lambda_{1})-\omega^{1}_{55}(\lambda_{5}-\lambda_{1})]=
\omega_{23}^{5}[\omega_{22}^{1}(\lambda_{3}-\lambda_{5})\\+\omega_{33}^{1}(\lambda_{3}-\lambda_{1})-\omega^{1}_{55}(\lambda_{5}-\lambda_{1})].
\end{array}
\end{equation}

Now, (\ref{3.47}) and (\ref{3.48}) is a homogeneous system of
equations in two variables $\omega_{32}^{5}$ and $\omega_{23}^{5}$
having either non trivial solution or trivial solution. If it has
trivial solution only, then we have $\omega_{32}^{5}=0$ and
$\omega_{23}^{5}=0$.

If it has non trivial solution also, then the determinant formed by
the coefficients of $\omega_{32}^{5}$ and $\omega_{23}^{5}$ in
(\ref{3.47}) and (\ref{3.48}) will be zero, i.e.,
\begin{center}
 $(\lambda_{3}-\lambda_{5})\omega_{22}^{1}+(\lambda_{5}-\lambda_{2})\omega_{33}^{1}+(\lambda_{2}-\lambda_{3})\omega_{55}^{1}=0$,
\end{center}
which gives
\begin{equation}\label{3.49}
\frac{\omega_{33}^{1}-\omega_{55}^{1}}{\lambda_{3}-\lambda_{5}}=\frac{\omega_{22}^{1}-\omega_{55}^{1}}{\lambda_{2}-\lambda_{5}}=k,
\end{equation}
where $k$ is constant.

Now, by using $\lambda_{5}=\frac{5H}{2}$ and (\ref{3.21}),
(\ref{3.49}) gives the expression
\begin{equation}\label{3.50}
2H\omega_{33}^{1}=kH(2\lambda_{3}-5H)+e_{1}(H).
\end{equation}

Again, differentiating (\ref{3.50}) along $e_{1}$ and using
(\ref{3.50}) and (\ref{3.51}), we get
\begin{equation}\label{3.52}
k(2\lambda_{3}-10H)e_{1}(H)-10H^{2}k^{2}(2\lambda_{3}-5H)-10H^{2}\lambda_{3}+25H^{3}=0
\end{equation}

After differentiating (\ref{3.52}) along $e_{1}$ and using
(\ref{3.50}) and (\ref{3.51}), we get
\begin{equation}\label{3.53}
\begin{array}{rcl}
k(8\lambda_{3}-45H)e_{1}^{2}(H)+He_{1}(H)[k^{2}(4\lambda_{3}^{2}-100\lambda_{3}H+225H^{2})-50\lambda_{3}H\\+125H^{2}]
-50H^{4}k(\lambda_{3}-5H)-5H^{2}k(4\lambda_{3}^{2}-25H^{2})(2k^{2}+1)=0
\end{array}
\end{equation}

Eliminating $e_{1}(H)$ from (\ref{3.53}), using (\ref{3.52}), we
obtain an algebraic equation in $\lambda_{3}$ and $H$ defined as
\begin{equation}\label{3.54}
 L(\lambda_{3},H)=0.
\end{equation}

If we differentiate (\ref{3.54}) along $e_{1}$ and using
(\ref{3.50}), a direct computation gives again an algebraic equation
in $\lambda_{3}$ and $H$ defined as
\begin{equation}\label{3.55}
 M(\lambda_{3},H)=0.
\end{equation}

We rewrite $L(\lambda_{3},H)$, $M(\lambda_{3},H)$ as polynomials
$L_{H}(\lambda_{3}), M_{H}(\lambda_{3})$ of $\lambda_{3}$ with
coefficients in the polynomial ring $R[H]$ over $\mathbb{R}$. Since
$L_{H}(\lambda_{3}) = M_{H}(\lambda_{3}) = 0$, $\lambda _3$ is a
common root of $L_H, M_H$, which implies that resultant $\Re(L_{H},
M_{H}) = 0$. It is obvious that $\Re(L_{H}, M_{H})$ is a polynomial
of $H$ with constant coefficients, therefore $H$ must be a constant
which is a contradiction to (\ref{3.3}). Hence
$\omega^{5}_{32}=\omega^{5}_{23}=0$. In similar way, we can prove
$\omega^{5}_{34}=\omega^{5}_{43}=\omega^{5}_{42}=\omega^{5}_{24}=0,$
which completes the proof of Lemma \ref{lemma3.3}. $\qed$
\medskip
\vspace{0.05in}\\

Now, using Lemma \ref{lemma3.3}, (\ref{3.5}) and (\ref{3.7}), we
have
\begin{equation}\label{3.56}
\omega_{5i}^{j}=\omega_{i5}^{j}=0 \quad i,j=2,3,4.
\end{equation}

Evaluating $g(R(X,Y)Z,W)$, using Lemma \ref{lemma3.1} and Lemma
\ref{lemma3.3}, Gauss equation
(\ref{2.5}) and (\ref{3.56}), we obtain the following:\\

 $\bullet\hspace{.2 cm}
 \mbox{For}\hspace{.2 cm} X=e_{2}, Y=e_{5}, Z=e_{2}, W= e_{5},$
\begin{equation}\label{3.60}
-\omega_{22}^{1}\omega_{55}^{1} = \frac{5H}{2}\lambda_{2}.
\end{equation}

 $\bullet\hspace{.2 cm}
 \mbox{For}\hspace{.2 cm} X=e_{3}, Y=e_{5}, Z=e_{3}, W= e_{5},$
\begin{equation}\label{3.61}
-\omega_{33}^{1}\omega_{55}^{1} =\frac{5H}{2} \lambda_{3}.
\end{equation}

 $\bullet\hspace{.2 cm}
 \mbox{For}\hspace{.2 cm} X=e_{4}, Y=e_{5}, Z=e_{4}, W= e_{5},$
\begin{equation}\label{3.62}
-\omega_{44}^{1}\omega_{55}^{1} = \frac{5H}{2}\lambda_{4}.
\end{equation}

Adding (\ref{3.60}), (\ref{3.61}), (\ref{3.62}) and using
(\ref{3.11}), we obtain
\begin{equation}\label{3.67}
\omega_{22}^{1}+\omega_{33}^{1}+ \omega_{44}^{1}=
-\frac{25H^{3}}{e_{1}(H)}.
\end{equation}

From (\ref{3.3}) and Lemma \ref{lemma3.1}, and the fact that
$[e_{i}\hspace{.1 cm}e_{1}](H)=0=
\nabla_{e_{i}}e_{1}(H)-\nabla_{e_{1}}e_{i}(H),$ for
 $i=2, 3, 4, 5$, we obtain
\begin{equation}\label{3.24}
 e_{i}e_{1}(H)= 0, \hspace{1 cm} i= 2, 3, 4, 5.
\end{equation}

 Using (\ref{3.67}) and (\ref{3.24}), we
  get
\begin{equation}\label{3.25}
 e_{i}( \omega_{22}^{1} + \omega_{33}^{1}+ \omega_{44}^{1})= 0, \hspace{1 cm} i= 2, 3,
 4, 5.
\end{equation}

Now, we need the following:

\begin{lem}\label{lemma3.2} Let $M^{5}$ be a $\delta(4)$-ideal
biconservative hypersurface of non constant mean curvature in
Euclidean space $\mathbb{E}^{6}$. Then,
 $e_{i}(\lambda_{j})=0$, for $i, j= 2,
3, 4$, and $i\neq j.$
\end{lem}
\noindent {\it Proof.} Operating with $e_{2}$ on both sides of
(\ref{3.22}), (\ref{3.11}) and using (\ref{3.6}), we find
\begin{equation}\label{3.26}
 (\lambda_{2}-\lambda_{4})^{2}\omega_{44}^{2}+(\lambda_{2}-\lambda_{3})^{2}\omega_{33}^{2}=0.
\end{equation}

Differentiating (\ref{3.26}) along $e_{1}$ and using (\ref{3.6}),
(\ref{3.19}) for $i,j=2,3$, we get
\begin{center}
$[-2(\lambda_{1}-\lambda_{2})(\lambda_{2}-\lambda_{4})\omega_{22}^{1}+(2\lambda_{1}+
 \lambda_{2}-3\lambda_{4})(\lambda_{2}-\lambda_{4})\omega_{44}^{1}]\omega_{44}^{2}
  +[-2(\lambda_{1}-\lambda_{2})(\lambda_{2}-\lambda_{3})\omega_{22}^{1}+(2\lambda_{1}+
 \lambda_{2}-3\lambda_{3})(\lambda_{2}-\lambda_{3})\omega_{33}^{1}]\omega_{33}^{2}
 =0.$
\end{center}
and using (\ref{3.26}) in the above equation, we get
\begin{equation}\label{3.27}
\begin{array}{rcl}
 [2(\lambda_{1}-\lambda_{2})(\lambda_{3}-\lambda_{4})\omega_{22}^{1}+(2\lambda_{1}+\lambda_{2}
 -3\lambda_{4})(\lambda_{2}-\lambda_{3})\omega_{44}^{1}
 \\-(2\lambda_{1}+\lambda_{2}-3\lambda_{3})(\lambda_{2}-\lambda_{4})\omega_{33}^{1}]\omega_{44}^{2}=0.
\end{array}
\end{equation}

Similarly, acting with $e_{1}$ and $e_{2}$ on (\ref{3.11}),
successively and using (\ref{3.6}), (\ref{3.20}) for $j=2$ and
$i=3,4$, (\ref{3.25}) and (\ref{3.26}), we obtain
\begin{equation}\label{3.28}
\begin{array}{rcl}
 [(\lambda_{4}-\lambda_{3})\omega_{22}^{1}+(\lambda_{3}-\lambda_{2})\omega_{44}^{1}
 + (\lambda_{2}-\lambda_{4})\omega_{33}^{1}]\omega_{44}^{2}=0.
\end{array}
\end{equation}

Equations (\ref{3.27}) and (\ref{3.28}) show that either
$\omega_{44}^{2}$, or the expression between square brackets, has to
vanish. We now prove that $\omega_{44}^{2}$ has to be zero. In fact,
if $\omega_{44}^{2}\neq 0$, then the expressions between square
brackets has to be zero:
\begin{equation}\label{3.29}
\begin{array}{rcl}
 2(\lambda_{1}-\lambda_{2})(\lambda_{3}-\lambda_{4})\omega_{22}^{1}+(2\lambda_{1}+\lambda_{2}-3\lambda_{4})(\lambda_{2}-\lambda_{3})\omega_{44}^{1}
 \\-(2\lambda_{1}+\lambda_{2}-3\lambda_{3})(\lambda_{2}-\lambda_{4})\omega_{33}^{1}=0.
\end{array}
\end{equation}
\begin{equation}\label{3.30}
\begin{array}{rcl}
 (\lambda_{4}-\lambda_{3})\omega_{22}^{1}+(\lambda_{3}-\lambda_{2})\omega_{44}^{1}
 + (\lambda_{2}-\lambda_{4})\omega_{33}^{1}=0.
\end{array}
\end{equation}

Eliminating $\omega_{22}^{1}$ from (\ref{3.29}) and (\ref{3.30}), we
get
\begin{equation}\label{3.31}
(\lambda_{2}-\lambda_{3})
(\lambda_{2}-\lambda_{4})(\omega_{33}^{1}-\omega_{44}^{1})=0,
\end{equation}
which shows that
\begin{equation}\label{3.32}
\omega_{33}^{1}-\omega_{44}^{1}=0,
\end{equation}
then using it to eliminate $\omega_{33}^{1}$, from (\ref{3.29}) and
(\ref{3.30}), we find
\begin{equation}\label{3.33}
\begin{array}{rcl}
 2(\lambda_{1}-\lambda_{2})(\lambda_{3}-\lambda_{4})\omega_{22}^{1}+
 [(2\lambda_{1}+\lambda_{2}-3\lambda_{4})(\lambda_{2}-\lambda_{3})\\+(2\lambda_{1}+\lambda_{2}-3\lambda_{3})(\lambda_{2}-\lambda_{4})]\omega_{44}^{1}=0.
\end{array}
\end{equation}
\begin{equation}\label{3.34}
 (\lambda_{4}-\lambda_{3})\omega_{22}^{1}+
 (\lambda_{3}+\lambda_{4}-2\lambda_{2})\omega_{44}^{1}=0.
\end{equation}

From (\ref{3.33}) and (\ref{3.34}), we obtain
\begin{equation}\label{3.35}
 (\lambda_{2}-\lambda_{3})(\lambda_{2}-\lambda_{4})=0,
\end{equation}
which contradicts the fact that principal curvatures are distinct.
Therefore, $\omega_{44}^{2}=0$, which gives $\omega_{33}^{2}=0$ in
view of (\ref{3.26}). Consequently,
$e_{2}(\lambda_{3})=e_{2}(\lambda_{4})=0$. \\

Similarly, we can prove that
$e_{3}(\lambda_{2})=e_{3}(\lambda_{4})=e_{4}(\lambda_{2})=e_{4}(\lambda_{3})=0$.
 which completes the
proof of Lemma \ref{lemma3.2}. \qed
\medskip
\vspace{0.05in}\\

Now, evaluating $g(R(X,Y)Z,W)$, using Lemma
\ref{lemma3.1}$\sim$\ref{lemma3.2}, Gauss equation
(\ref{2.5}) and (\ref{3.56}), we obtain the following:\\

  $\bullet \hspace{.2 cm}
 \mbox{For}\hspace{.2 cm} X=e_{2}, Y=e_{3}, Z=e_{2}, W= e_{3}$
\begin{equation}\label{3.57}
\begin{array}{rcl}
-\omega_{22}^{1}\omega_{33}^{1}+\omega_{32}^{4}\omega_{23}^{4}
-\omega_{34}^{2}\omega_{43}^{2}-\omega_{42}^{3}\omega_{24}^{3}=
\lambda_{2} \lambda_{3}.\end{array}
\end{equation}

  $\bullet\hspace{.2 cm}
 \mbox{For}\hspace{.2 cm} X=e_{2}, Y=e_{4}, Z=e_{2}, W= e_{4}$
\begin{equation}\label{3.58}
\begin{array}{rcl}
-\omega_{22}^{1}\omega_{44}^{1}
+\omega_{42}^{3}\omega_{24}^{3}-\omega_{32}^{4}\omega_{23}^{4}
-\omega_{34}^{2}\omega_{43}^{2}= \lambda_{2} \lambda_{4}.
\end{array}
\end{equation}

  $\bullet\hspace{.2 cm}
 \mbox{For}\hspace{.2 cm} X=e_{3}, Y=e_{4}, Z=e_{3}, W= e_{4},$
\begin{equation}\label{3.59}
\begin{array}{rcl}
-\omega_{33}^{1}\omega_{44}^{1}
+\omega_{34}^{2}\omega_{43}^{2}-\omega_{32}^{4}\omega_{23}^{4}
-\omega_{42}^{3}\omega_{24}^{3}= \lambda_{3} \lambda_{4}.\end{array}
\end{equation}

Now, using (\ref{3.7}) and (\ref{3.22}), we have
\begin{equation}\label{3.63}
 \lambda_{2} \lambda_{3}+ \lambda_{3} \lambda_{4}+ \lambda_{4}
 \lambda_{2}=\frac{25H^{2}}{4}+\frac{\rho}{2}.
\end{equation}

Also, from (\ref{3.5}) and (\ref{3.7}), we have
\begin{equation}\label{3.64}
(\lambda_{2}-\lambda_{3})\omega_{42}^{3}=(\lambda_{4}-\lambda_{3})\omega_{24}^{3}=(\lambda_{2}-\lambda_{4})\omega_{32}^{4},
\end{equation}
which gives
\begin{equation}\label{3.65}
\omega_{42}^{3}\omega_{24}^{3}+\omega_{43}^{2}\omega_{34}^{2}+\omega_{32}^{4}\omega_{23}^{4}=0.
\end{equation}

Adding (\ref{3.56}),(\ref{3.57}),(\ref{3.58}) and using
(\ref{3.63}), (\ref{3.65}), we get
\begin{equation}\label{3.66}
\omega_{22}^{1}\omega_{33}^{1}+
\omega_{44}^{1}\omega_{33}^{1}+\omega_{22}^{1}\omega_{44}^{1}=
-\frac{25H^{2}}{4}-\frac{\rho}{2}.
\end{equation}

By squaring (\ref{3.60}), (\ref{3.61}), (\ref{3.62}) and then
adding, we have the expression
\begin{equation}\label{3.68}
(\omega_{55}^{1})^{2}[(\omega_{22}^{1})^{2}+(\omega_{33}^{1})^{2}+
(\omega_{44}^{1})^{2}]=
\frac{25H^{2}}{4}(\lambda_{2}^{2}+\lambda_{3}^{2}+\lambda_{4}^{2}).
\end{equation}

By direct calculation, using (\ref{3.21}), (\ref{3.22})
(\ref{3.66}), (\ref{3.67}) and (\ref{3.68}), we find
\begin{equation}\label{3.69}
e_{1}^{2}(H)=-25H^{4},
\end{equation}
which is possible only when $H=0.$\\

\textbf{Case B.} $\lambda_{2}=\lambda_{3}$ and $\lambda_{3},
\lambda_{4}, \lambda_{5}$ are distinct.

In this case, from (\ref{3.11}) and (\ref{3.22}) we have
\begin{equation}\label{3.70}
2\lambda_{3}+\lambda_{4}=5H,
\end{equation}
\begin{equation}\label{3.71}
2\lambda_{3}^{2}+\lambda_{4}^{2}=\frac{25H^{2}}{2}-\rho.
\end{equation}

Using (\ref{3.66}) and (\ref{3.67}), we find
\begin{equation}\label{3.72}
6\lambda_{3}^{2}-20H\lambda_{3}+\frac{25H^{2}}{2}+\rho=0,
\end{equation}

Differentiating (\ref{3.68}), (\ref{3.69}) along $e_{1}$ and using
(\ref{3.6}) and (\ref{3.66}), we obtain
\begin{equation}\label{3.73}
e_{1}(\lambda_{3})=\frac{20\lambda_{3}-25H}{12\lambda_{3}-20H}e_{1}(H),
\end{equation}
\begin{equation}\label{3.74}
\omega^{1}_{33}=\frac{20\lambda_{3}-25H}{(2\lambda_{3}+5H)(6\lambda_{3}-10H)}e_{1}(H),
\end{equation}

Differentiating (\ref{3.74}) along $e_{1}$ and using (\ref{3.74})
and (\ref{3.6}), we obtain
\begin{equation}\label{3.75}
\begin{array}{rcl}
e_{1}^{2}(H)[89750\lambda_{3}H^{4} -16400\lambda_{3}^{2}H^{3}
-23190\lambda_{3}^{3}H^{2} + 5016\lambda_{3}^{4}H +
1800\lambda_{3}^{5}\\ -55625H^{5}]= 2H^{2}(-125000\lambda_{3}H^{6} +
9375\lambda_{3}^{2}H^{5} + 53750\lambda_{3}^{3}H^{4}
\\-10250\lambda_{3}^{4}H^{3} -8100\lambda_{3}^{5}H^{2} +
1080\lambda_{3}^{6}H + 432\lambda_{3}^{7} + 78125H^{7}) ,
\end{array}
\end{equation}

Differentiating (\ref{3.74}) along $e_{1}$ and using (\ref{3.73}),
we obtain
\begin{equation}\label{3.76}
\begin{array}{rcl}
e_{1}^{2}(H)P(H,\lambda_{3})=Q(H,\lambda_{3}).\end{array}
\end{equation}

Eliminating $e_{1}^{2}(H)$ from (\ref{3.75}) and (\ref{3.76}), we
get
\begin{equation}\label{3.77}
F(\lambda_{3}, H)=0,
\end{equation}
which is the polynomial equation in $H$ and $\lambda_{3}$. It can be
easily seen that if we differentiate (\ref{3.77}) along $e_{1}$ and
using (\ref{3.73}), we find another polynomial equation in $H$ and
$\lambda_{3}$ defined as
\begin{equation}\label{3.78}
G(\lambda_{3}, H)=0.
\end{equation}

We rewrite $F(\lambda_{3},H)$, $G(\lambda_{3},H)$ as polynomials
$F_{H}(\lambda_{3}), G_{H}(\lambda_{3})$ of $\lambda_{3}$ with
coefficients in the polynomial ring $R[H]$ over $\mathbb{R}$. Since
$F_{H}(\lambda_{3}) = G_{H}(\lambda_{3}) = 0$, $\lambda _3$ is a
common root of $F_H, G_H$, which implies that resultant $\Re(F_{H},
G_{H}) = 0$. It is obvious that $\Re(F_{H}, G_{H})$ is a polynomial
of $H$ with constant coefficients, therefore $H$ must be a
constant.\\

\textbf{Case C.} $\lambda_{2}= \lambda_{3} = \lambda_{4}\neq
\lambda_{5}$

In this case, using (\ref{3.11}) and (\ref{3.22}), we find that
$\rho =
-\frac{25H^{2}}{2}$, which implies that $H$ is a constant.\\

Combining all above cases \textbf{A}, \textbf{B} and \textbf{C}, we
can conclude Theorem \ref{main}.\\

 \noindent
\textbf{Acknowledgement}. \small \emph{The authors acknowledge
useful discussions and suggestions with Dr. Ram Shankar Gupta. Also,
the first author acknowledge the support provided by NISER,
Bhubaneswar. The second author was supported by Grant $\# E.037$
from the Research Committee of the University of Patras (Programme
K. Karatheodori).}


\bibliography{xbib}


Author's address:\\
\\
\textbf{Deepika}\\
S N Bose National Centre for Basic Sciences,\\
Sector III, Kolkata- 700106, India.\\
\textbf{Email:} sdeep2007@gmail.com\\
\\
\textbf{Andreas Arvanitoyeorgoes}\\
University of Patras,\\
Department of Mathematics,\\
 GR-26500 Patras, Greece.\\
 \textbf{Email:} arvanito@math.upatras.gr\\

\end{document}